# Structural shrinkage of nonparametric spectral estimators for multivariate time series

Hilmar Böhm

Rainer von Sachs[*]

**Abstract:** In this paper we investigate the performance of periodogram based estimators of the spectral density matrix of possibly high-dimensional time series. We suggest and study shrinkage as a remedy against numerical instabilities due to deteriorating condition numbers of (kernel) smoothed periodogram matrices. Moreover, shrinking the empirical eigenvalues in the frequency domain towards one another also improves at the same time the Mean Squared Error (MSE) of these widely used nonparametric spectral estimators. Compared to some existing time domain approaches, restricted to i.i.d. data, in the frequency domain it is necessary to take the size of the smoothing span as "effective or local sample size" into account. While Böhm and von Sachs (2007) proposes a multiple of the identity matrix as optimal shrinkage target in the absence of knowledge about the multidimensional structure of the data, here we consider "structural" shrinkage. We assume that the spectral structure of the data is induced by underlying factors. However, in contrast to actual factor modelling suffering from the need to choose the number of factors, we suggest a model-free approach. Our final estimator is the asymptotically MSE-optimal linear combination of the smoothed periodogram and the parametric estimator based on an underfitting (and hence deliberately misspecified) factor model. We complete our theoretical considerations by some extensive simulation studies. In the situation of data generated from a higher-order factor model, we compare all four types of involved estimators (including the one of Böhm and von Sachs (2007)).

## 1. Introduction

Spectral analysis of multivariate time series is known to be a useful tool to analyse not only serial but also cross-correlations of dynamic data of possibly high dimension (Shumway and Stoffer, 2000). In the absence of some possibly restrictive parametric assumptions on the dynamics of the time series (such as vector autoregressive - moving average of finite order), the standard nonparametric approach of smoothing the periodogram matrix over frequency usually shares well-established and generally even for moderate sample sizes satisfactory properties such as approximate unbiasedness, approximate uncorrelatedness over different frequencies and the usual variance-bias trade off known from classical nonparametric theory (Brillinger (1975)). What is less known and explored, however, and highly relevant for more and more frequently met situations of large dimensionality of the time series, is the deterioration of the condition number of the resulting nonparametric estimator (smoothed periodogram matrix). It is known that a high condition number of such a matrix, i.e. the ratio $l_{\max}/l_{\min}$ of its largest to its smallest eigenvalue, leads to numerical instabilities, in particular when the (estimated) spectral density matrix is used subsequently in sensitive functionals such as its inverse or its determinant. A prominent example for the latter ones is the use of the Kullback-Leibler discrimination information (Kullback and Leibler, 1952), as a measure of disparity between several estimated multivariate spectra (as in Kakizawa, Shumway and Taniguchi (1998), e.g.), to be used in classification of multivariate time series.

---

[*]We acknowledge financial support from the IAP research network grant P 5/24 of the Belgian government (Belgian Science Policy) as well as from the contract "Projet d'Actions de Recherche Concertées" n° 07/12-002 of the "Communauté française de Belgique", granted by the "Académie universitaire Louvain".

In many fields of application, including economic panel data (Bai and Ng, 2002; Forni, Hallin, Lippi and Reichlin, 2000), but also genetic engineering or neuropsychology, the dimension of the data can come close to the sample size, making the smoothed periodogram become close to a singular matrix, in particular.

In this paper we suggest a remedy to improve upon the smoothed periodogram as an estimator for the multivariate spectrum using regularization, i.e. shrinkage, techniques. It is known from the statistical literature on estimation in i.i.d. data situations (Haff, 1977, 1979, 1980), that shrinkage helps to correct the following effect: the dispersion of the sample eigenvalues can be tremendously larger than the dispersion of the population eigenvalues of the spectrum as the large eigenvalues are biased upwards, the small ones downwards (Jolliffe, 2002). Thus, the quality of an estimator of a high-dimensional target can be improved, by shrinking the eigenvalues towards one another, not only numerically, but even on the level of the widely used criterion of mean square error (Beran and Dümbgen, 1998; Ledoit and Wolf, 2004).

We note that this technique is more than just standardizing each dimension of the time series - which would improve the condition number in case of minimally coherent data, but not so with potentially highly cross-correlated data (the interdependence over dimension being responsable for the afore-mentioned dispersion effect).

Compared to existing work on shrinkage in the time domain, we show that in the frequency domain it is necessary to take the size of the smoothing span $m$ as "effective or local sample size" into account. We note that simply choosing the smoothing span of the smoothed periodogram sufficiently large is no reasonable solution to the problem: depending on the roughness of the true spectral density to be estimated, this might result into important oversmoothing.

For reasons of notational simplicity, in this work, we consider as simplest smoothing method the averaged periodogram, that is a symmetric kernel smoother of finite support ("boxcar") with equal weights for each periodogram ordinate within the smoothing span. One can easily check that all the results of our paper carry directly over to the more frequently used kernels in the smooth-ing literature.

Our proposed shrinkage estimator is, pointwise at frequency $\omega \in (0, 2\pi]$, a convex combination of the averaged periodogram $\hat{f}_T^0(\omega)$ with some shrinkage target $\hat{f}_T^1(\omega)$ in the frequency domain. I.e., our estimators are of the form $\hat{f}_T(\omega) := r_T(\omega)\hat{f}_T^1(\omega) + (1 - r_T(\omega))\hat{f}_T^0(\omega)$, where in order to reduce the dispersion of the eigenvalues of $\hat{f}_T^0(\omega)$, the factor $r_T$ is chosen such that the sample eigenvalues are shrunk towards each other linearly. The most direct target to use would be (a multiple of) the identity matrix, i.e. $\hat{f}_T^1(\omega) = \mu(\omega)$ Id. This set-up has been treated by the authors in a companion paper (Böhm and von Sachs, 2007), where they determine the optimal amount of shrinkage by a data driven approach in a framework of an asymptotically with sample size growing dimension.

Obviously, the technique of shrinkage has a certain relationship to ridge regression. In fact, a linear combination of a sample covariance and the identity matrix has been used as original motivation for ridge regression (Hoerl and Kennard, 1970). However, shrinkage of empirical covariance matrices or spectral density estimators is not the same thing as ridge estimation. In the first approach, the eigenvalues of the matrices under consideration are shrunken towards each other, and hence their dispersion is reduced. Constructing a ridge, all eigenvalues are moved away from zero (either by the same amount for ordinary ridge or by some individual constant for generalized ridge regression) in order to regularize the estimator. Recent theoretical work of Bickel and Levina (2007), Rothman, Levina and Zhu (2008), e.g., using a lot more refined techniques such as Lasso, or thresholding for regularization of large-dimensional covariance matrices, are in this latter spirit.

Using the identity matrix as a shrinkage target is reasonable if there is little or no knowledge about the underlying multidimensional structure of the data. In this case, a shrinkage target should be used that imposes the least possible amount of structure and which, at the same time, has the best of all possible condition numbers. In many settings, however, it is reasonable to assume that the covariance or spectral structure of the data is induced by underlying, known or hidden, factors. The general idea underlying factor models is that $p$ observed random variables can be expressed, except for an error term, as linear functions of $q < p$ random factors. For instance, in econometrics, markets are usually assumed to be driven by underlying global variables such as interest rate, employment rate or gross national product. The models reach from simple one-factor models, as in Sharpe (1963), to sophisticated approaches that use multiple global and industry specific factors that may be intercorrelated, as, e.g., in Forni et al. (2000).

A disadvantage of factor models is that, usually, the number of factors is a parameter that must be either specified *a priori* or chosen by somewhat sophisticated data-driven procedures akin to model selection. Research on how to propose a generally satisfying criterion is still going on (Bai and Ng, 2002; Hallin and Liška, 2007), and it would be interesting to avoid this problem while taking advantage of the structure imposed by a factor model to be a remedy to the curse of dimensionality.

We have developed a hybrid approach to circumvent the dilemma of model choice and still retain the advantages of factor analysis. We combine a nonparametric estimator, in our case the averaged periodogram $\hat{f}_T^0(\omega)$, with a parametric estimator $\hat{f}_T^1(\omega)$ of the spectral matrix. The latter is our new shrinkage target. It is given by fitting a one-factor model to the data. However, we do not assume that this model is true; rather, we believe that the data follow a more complicated structure. This may be a $q$-factor model ($q > 1$), a model driven by different layers of factors, or the model may be completely unknown. By combining a shrinkage target, which is actually underfitted, with a nonparametric estimator of the spectrum, we circumvent the problem of model choice. In a data driven approach, weights are chosen such that the new, hybrid estimator is the asymptotically optimal linear combination of two conventional estimators. The first component, the averaged periodogram, is asymptotically unbiased but has high variance. The second component is biased due to misspecification but, by imposing structure, has low variance.

We note that, instead of choosing a one-factor model as our shrinkage target, we might as well opt for something more complicated, e.g. a $q$ factor model with $q > 1$. The only prerequisite for doing this is having background knowledge that the underlying structure is more complicated than the shrinkage target, e.g. a $\tilde{q}$ factor model, $\tilde{q} > q$. The theory we will give in section 3.1 can easily be adapted to such a case.

To the best of the authors' knowledge, there is no literature on shrinkage to a factor model in time series analysis. In the literature on finance, an approach to shrink to a factor model has been developed in the context of portfolio selection (Ledoit and Wolf, 2003) under iid assumptions on the data. However, the idea of shrinking a nonparametric fit towards a parametric estimator drives quite generally a variety of existing approaches, among which one finds the work of Daniels and Cressie (2001), and to some extent, of Botts and Daniels (2006), in the context of Bayesian covariance and spectral estimation, respectively.

The remainder of this paper is organized as follows: in the next section, we will develop the theoretical background for data driven shrinkage to a 'market' one-factor model, where the term 'market' is just a wildcard term that does not necessarily mean that we are in an economic context. We will first give the basic assumptions and definitions in the following subsection. In section 3.2, we will introduce the shrinkage target, which is a one-factor model. The model assumptions are that the $p$ dimensional process is driven by a dynamic, hidden or known, underlying process with spectral density $f_0(\omega)$. We will fit this model to the data; however, at the same time we assume that the model be misspecified. The philosophy behind this is that the model is just a parsimonious tool of describing the data. In sections 3.4 and 3.5 we derive the MSE-optimal solution for the shrinkage intensity which is a function of the true spectral density $f(\omega)$. In section 3.6, we will examine the asymptotic behaviour of the MSE-optimal shrinkage intensity $r_T(\omega)$, which will help us to develop a data driven estimator in section 3.7. Comprehensive Monte Carlo studies will show the usefulness of our estimator in section 4. We note that most proofs are relegated to an appendix section.

## 2. Multivariate spectral analysis

We assume that we observe a realization $(X_t)_{t=1}^T$ of a $p$-dimensional real-valued, centered stationary Gaussian time series $(X_t)$. We aim at estimating the $p \times p$ spectral density matrix function at frequency $\omega \in (0, 2\pi]$

$$f(\omega) = \frac{1}{2\pi T} \sum_{u \in \mathbb{Z}} \text{Cov}(X_t, X_{t+u}) \exp(-\iota \omega u), \quad \omega \in (0, 2\pi] \tag{1}$$

where $\iota = \sqrt{-1}$. The most common nonparametric estimators of (1) are based on the *periodogram*. If we denote by

$$d_T(\omega) = \frac{1}{\sqrt{2\pi T}} \sum_{t=1}^T X_t \exp(-\iota \omega t), \quad \omega \in (0, 2\pi] \tag{2}$$

the vector-valued *discrete Fourier transform* of the realization $(X_t)_{t=1}^T$, then the $p \times p$ periodogram matrix is defined as

$$I_T(\omega) := d_T(\omega) d_T^*(\omega) \tag{3}$$

where $^*$ means conjugate complex transpose. Furthermore, we will denote conjugate complex (for a scalar value) by overline. The periodogram is not a consistent estimator of the spectrum (1), but it is asymptotically unbiased. Moreover, for $p > 1$, the periodogram is a singular matrix: if $d_T(\omega) = (d_1(\omega), \ldots, d_p(\omega))'$, then (3) can be expressed as

$$I_T(\omega) = \left( \overline{d_1(\omega)} \begin{pmatrix} d_1(\omega) \\ \vdots \\ d_p(\omega) \end{pmatrix} \quad \ldots \quad \overline{d_p(\omega)} \begin{pmatrix} d_1(\omega) \\ \vdots \\ d_p(\omega) \end{pmatrix} \right) \tag{4}$$

and thus has almost surely rank 1. If the periodogram is smoothed over frequency, the estimators derived this way are consistent under a classical asymptotical framework. We will restrict ourselves to the simplest form of smoothing, the *averaged periodogram* with smoothing span $m_T$, where the conditions $m_T/T \to 0$ and $m_T \to \infty$ as $T \to \infty$ guarantee consistency and asymptotic unbiasedness:

$$\hat{f}_T^0(\omega) := \frac{1}{m_T} \sum_{k=-(m_T-1)/2}^{(m_T-1)/2} I_T(\omega + \omega_k), \tag{5}$$

where $\omega_k$ denotes the Fourier frequency $2\pi k/T$.

## 3. Theoretical framework

### 3.1. Setup and assumptions

Our aim is to estimate the spectrum $f(\omega)$ of a $p$-variate Gaussian time series. We assume that we have realizations

$$(X_{it})_{t \in \{1,\ldots,T\}} = X_{i1}, \ldots, X_{iT}, \qquad i = 1, \ldots, p$$

Moreover, we assume that we have realizations from another, one dimensional time series

$$(X_{0t})_{t \in \{1,\ldots,T\}} = X_{01}, \ldots, X_{0T}$$

to which we refer as the *market* or *exogenous* time series. The market time series is thought to be a process that has a certain explanatory value for the other time series $(X_{it}), i = 1, \ldots, p$. One possible choice is to use the average over dimension in the time domain of the $(X_{it})_{i=1,\ldots,p}$,

$$X_{0t} = \frac{1}{p} \sum_{i=1}^{p} X_{it}$$

However, we make no special assumptions on the market time series. It would as well be possible to choose an external variable or the first principal component of the data.

We make the following assumptions:

*Assumption* 3.1. All our time series, including the market time series, are centered

$$\operatorname{E} X_{it} = 0 \qquad i = 0, \ldots, p$$

and stationary.

In this paper, purely for reasons of simplifying the presentation, we do not present our estimation results in terms of the spectrum directly, but rather choosing the expected periodogram

$$f_T^0(\omega) := \operatorname{E} \hat{f}_T^0(\omega) \tag{6}$$

as estimation target. This is possibly without loss of generality because the expected periodogram $f_T^0(\omega)$ approaches the true spectrum $f(\omega)$ with a rate of convergence suffiently fast to enable us to carry over our proofs immediately to estimate $f(\omega)$. In order to do so we make the following assumption:

*Assumption* 3.2. If

$$\gamma_{ij}(h) = \operatorname{E} X_{it} X_{j,t+h}, \qquad i,j = 0, \ldots, p, \quad h \in \mathbb{Z}$$

denotes the autocovariance function, then

$$\sum_{h=-\infty}^{\infty} |h|\, |\gamma_{ij}(h)| < \infty \qquad \forall\, i,j = 0, \ldots, p$$

Then, we have the following well-known result from Brillinger (1975) or Shumway and Stoffer (2000):

*Lemma* 3.1. Under assumption 3.2, $f(\omega)$ has (elementwise and for real- and imaginary parts separately) continuous derivatives of order one, and hence

$$f_T^0(\omega) - f(\omega) = \mathrm{O}\left(m_T/T\right).$$

The enhanced estimator we want to construct is gained by linearly combining a standard nonparametric estimator, in our case the averaged periodogram, with a shrinkage target. The latter is gained by fitting a one-factor model to the data, where the time series $X_{0t}$ is assumed to be the underlying factor.

We assume the dimension $p$ to be fixed while still $T \to \infty$. We denote the $i$th component of the discrete Fourier transform of the data at frequency $\omega$ as $d_i(\omega)$.

We furthermore make the following notational convention: whenever we use vector- or matrix valued terms, we will mean the respective $p$-dimensional vector or the $p \times p$ matrix unless we explicitly state otherwise. Thus, $f(\omega), f_T^0(\omega)$ and $\hat{f}_T^0(\omega)$ refer to the spectrum, expected averaged periodogram and averaged periodogram, respectively, of the time series $(X_{it})_{i=1,\ldots,p}$. We will also refer to the $p$-dimensional vector of the time series at time $t$ as $X_t$. However, when we look at components, we will use the index value zero to refer to the market time series. E.g., we refer to the cross-spectrum between the market series and the first component of $X_t$ as $f_{01}(\omega)$.

### 3.2. One-factor model

The shrinkage target is given by fitting a one-factor model to the data $(X_{it}), i = 1,\ldots,p$, which we will define in this section. We will use a different notation for the random variables to emphasize that this model is not assumed to hold true for the data $X_{it}$. Rather, we use the model as a parsimonious tool to approximate the spectral structure of the process.

Let us assume that we have a univariate *exogenous* time series $\dot{X}_{0t}, t = 1,\ldots,T$ with spectrum $\dot{f}_0(\omega)$. When we speak of exogenous, we mean that this data $\dot{X}_{0t}$ can be used as a factor time series that has some explicative value for the data in the sense of the following model:

$$\dot{X}_{it} = \beta_i \dot{X}_{0t} + \epsilon_{it} \quad i = 1,\ldots,p \tag{7}$$

The weights $\beta_i \in \mathbb{R}$ are non-random. The idiosyncratic components $\epsilon_{it}$ are assumed to be normally distributed and independent over time and dimension, and independent of $(\dot{X}_{0t})$:

$$\epsilon_{it} \sim \mathcal{N}\left(0, 2\pi(\sigma_i^\epsilon)^2\right) \tag{8}$$

In this simple factor model, all serial correlation in the data $\dot{X}_{it}$ originates from serial correlation in the exogenous time series $\dot{X}_{0t}$. The serial correlation of the exogenous time series is determined by its spectrum $\dot{f}_0(\omega)$.

The fact that the idiosyncratic components are uncorrelated over time and dimension is important, as in either other case, it would be impossible to identify the model under classical asymptotics (Forni et al., 2000). Together with the independence between the idiosyncratic components and the exogenous time series, this has two more advantages: first, it will allow us to use linear regression to estimate the $\beta_i$ and the $(\sigma_i^\epsilon)^2$. Second, this model implies, simply by linearity, the following relationship for the DFTs of the data:

*Lemma 3.2.*
$$\dot{d}_i(\omega) = \beta_i \dot{d}_0(\omega) + \dot{d}_i^\epsilon(\omega) \tag{9}$$

where $\dot{d}_i^\epsilon(\omega)$ is the DFT of the idiosyncratic components. Furthermore,

$$\dot{d}_i^\epsilon(\omega) \sim \mathcal{N}^C\left(0, (\sigma_i^\epsilon)^2\right) \quad \forall\, \omega \tag{10}$$

We see from (10) that the variance in the idiosyncratic components is independent of frequency. Furthermore, the weights $\beta = (\beta_1, \ldots, \beta_p)'$ are independent of the frequency, too, due to (7). This means that the spectrum under the above specified one-factor model (7) is

$$\tilde{f}^1(\omega) = \beta\beta' \dot{f}_0(\omega) + \Delta \tag{11}$$

where

$$\Delta = \begin{pmatrix} (\sigma_1^\epsilon)^2 & \cdots & 0 \\ \vdots & \ddots & \vdots \\ 0 & \cdots & (\sigma_p^\epsilon)^2 \end{pmatrix} \tag{12}$$

When it comes to estimation of the one-factor model, we will as afore-mentioned identify the spectrum with the expected averaged periodogram. Thus, instead of using the model (11), we will use the slightly modified model

$$\dot{f}^1(\omega) = \beta\beta' \dot{f}_0^0(\omega) + \Delta, \tag{13}$$

where $\dot{f}_0^0(\omega)$ means the expected averaged periodogram of the factor time series $\dot{X}_{0t}$.

### 3.3. Estimation of one-factor model

The model (13) is assumed *not* to hold true. However, even under these circumstances, it is possible to fit this model to the time series $X_{it}$ by choosing weights $\beta_i$ such that the $L_2$ distance between $f^0(\omega)$ and $\beta\beta' f_0(\omega)$ becomes minimal.

We will refer to this minimum $L_2$ distance spectral density under the one-factor model as to $f^1(\omega)$.

The fact that both weights $\beta_i$ and idiosyncratic variances $(\sigma_i^\epsilon)^2$ are independent of lag and frequency, respectively, enables us to estimate these parameters with standard methods. We use linear regression to obtain the following estimators $b_i$ for the $\beta_i$s

$$b_i = \frac{\sum_{t=1}^T (X_{0t} X_{it})}{\sum_{t=1}^T (X_{0t})^2}, \tag{14}$$

which is just the standard estimator of the slope in linear regression.

Next, we need to estimate the variances $(\sigma_i^\epsilon)^2$ of the idiosyncratic components. The standard way to do this is again to use the time domain estimator of the residual variance, which we normalize by $1/2\pi$:

$$\widehat{(\sigma_i^\epsilon)^2} = \frac{1}{2\pi} \sum_{t=1}^{T} \frac{(X_{0t} - b_i X_{it})^2}{T} \qquad (15)$$

Furthermore, both estimators have the convenient property of being consistent, and the stochastic rate of convergence is in both cases $1/\sqrt{T}$ (Sachs and Hedderich, 2006):

$$b_i = \beta_i + \mathrm{O}_p\left(\frac{1}{\sqrt{T}}\right) \qquad (16)$$

and

$$\widehat{(\sigma_i^\epsilon)^2} = (\sigma_i^\epsilon)^2 + \mathrm{O}_p\left(\frac{1}{\sqrt{T}}\right) \qquad (17)$$

Plugging the estimators from (16) and (17) and the averaged periodogram of $X_{0t}$,

$$\hat{f}_0^0(\omega) = \frac{1}{m_T} \sum_{k=-(m_T-1)/2}^{(m_T-1)/2} I_{00}(\omega + \omega_k) ,$$

into the definition of the one-factor model (13), we obtain an estimator of the multivariate spectrum that is based on a one-factor model:

$$\hat{f}_T^1(\omega) = bb' \hat{f}_0^0(\omega) + D , \qquad (18)$$

where

$$D = \mathrm{diag}\left(\widehat{(\sigma_1^\epsilon)^2} \ldots \widehat{(\sigma_p^\epsilon)^2}\right) .$$

This estimator is our shrinkage target. By construction of the model, with equations (14) and (15), we observe that on the diagonal $\hat{f}_T^1(\omega) = \hat{f}_T^0(\omega)$.

### 3.4. Optimal shrinkage intensity

Our aim is to improve upon the averaged periodogram by shrinking to a target matrix function that is more regular, at the price of possibly having larger bias. Here, we make the assumption that a one-factor model is not far too crude an approximation. We do, however, not believe that the underlying structure is totally explained; we even make the opposite assumption, namely that the model is misspecified:

*Assumption* 3.3. There exists a $\delta > 0$ such that, uniformly over all frequencies $\omega \in [0, 2\pi]$ and all $i, j = 1, \ldots, p$, we have

$$\left|f_{ij}^1(\omega) - f_{ij}^0(\omega)\right| \geq \delta \qquad (19)$$

Assumption 3.3 is made for technical reasons: the estimator of the shrinkage intensity which we are going to derive will have an estimator of the difference $f^1_{ij}(\omega) - f^0_{ij}(\omega)$ in the denominator. Because of this, assumption 3.3 is needed to avoid problems of identifiability.

We search for a linear combination

$$\hat{f}^+(\omega) = \zeta_T(\omega)\hat{f}^1_T(\omega) + (1 - \zeta_T(\omega))\hat{f}^0_T(\omega)$$

where $\zeta_T(\omega)$ is a data driven estimator of an optimal, oracle shrinkage intensity $\zeta^*_T(\omega)$ that is the solution of the minimization problem

$$\mathrm{E}\left\|\hat{f}^+_T(\omega) - f^0_T(\omega)\right\|^2 = \min!\ , \qquad (20)$$

that is,

$$\zeta^*_T(\omega) = \arg\min_{z_T(\omega)} \mathrm{E}\left\|z_T(\omega)\hat{f}^1_T(\omega) + (1 - z_T(\omega))\hat{f}^0_T(\omega) - f^0_T(\omega)\right\|^2.$$

We will proceed in three steps:

First, in subsection 3.5, we will derive the optimal, *oracle* shrinkage intensity $\zeta^*_T(\omega)$ which depends on background knowledge of the underlying process.

Second, in subsection 3.6 we will derive the asymptotic behaviour of the oracle shrinkage intensity. We will see that the necessity to shrink vanishes asymptotically. This is because the averaged periodogram is a consistent estimator whereas the shrinkage target is misspecified due to assumption 3.3. As a consequence, the data driven estimator of $f^0_T(\omega)$ will asymptotically have the same behaviour as the averaged periodogram, as the data driven estimator of the shrinkage intensity will converge to zero. Finally, we will construct a data driven estimator in subsection 3.7.

### *3.5. Oracle shrinkage intensity*

We will derive the *oracle* shrinkage intensity by solving the minimization problem given in formula (20). This can simply be done by differentiation. Let $z \in [0, 1]$ denote a shrinkage intensity. The risk $R(z)$ associated with $z$ is derived in Appendix A.1:

$$\begin{aligned} R(z) &= \mathrm{E}\left\|z\hat{f}^1_T(\omega) + (1-z)\hat{f}^0_T(\omega) - f^0_T(\omega)\right\|^2 \\ &= \sum_{i,j=1}^{p}\Big(z^2\,\mathrm{Var}\,\hat{f}^1_{ij}(\omega) + (1-z)^2\,\mathrm{Var}\,\hat{f}^0_{ij}(\omega) \\ &\qquad + 2z(1-z)\Re\left(\mathrm{Cov}\left(\hat{f}^1_{ij}(\omega), \hat{f}^0_{ij}(\omega)\right)\right) \\ &\qquad + z^2\left|f^1_{ij}(\omega) - f^0_{ij}(\omega)\right|^2\Big) \end{aligned} \qquad (21)$$

where we have used that $\mathrm{E}\,\hat{f}^0_T(\omega) = f^0_T(\omega)$
and, according to (13), $\mathrm{E}\,\hat{f}^1_T(\omega) = f^1_T(\omega)$.

The first derivative of R(z) with respect to $z$ is:

$$R'(z) = 2 \sum_{i,j=1}^{p} \left( z \operatorname{Var} \hat{f}_{ij}^1(\omega) - (1-z) \operatorname{Var} \hat{f}_{ij}^0(\omega) \right.$$
$$\left. + (1-2z)\Re\left(\operatorname{Cov}\left(\hat{f}_{ij}^1(\omega), \hat{f}_{ij}^0(\omega)\right)\right) + z\left|f_{ij}^1(\omega) - f_{ij}^0(\omega)\right|^2\right)$$

Moreover, the second derivative is

$$R''(z) = 2 \sum_{i,j=1}^{p} \left( \operatorname{Var}(\hat{f}_{ij}^1(\omega) - \hat{f}_{ij}^0(\omega)) + \left|f_{ij}^1(\omega) - f_{ij}^0(\omega)\right|^2 \right)$$
$$> 0 \qquad (22)$$

where we use that $\hat{f}_T^1(\omega)$ and $\hat{f}_T^0(\omega)$ are hermitian, so that the imaginary parts sum to zero.

Thus, we know that any local extremum will be a minimum. Setting the first derivative equal to zero, we obtain the following theorem.

*Theorem* 3.3. The optimal shrinkage intensity is given by

$$\zeta_T^*(\omega) = \frac{\sum_{i,j=1}^{p} \left( \operatorname{Var} \hat{f}_{ij}^0(\omega) - 2\Re \operatorname{Cov}\left(\hat{f}_{ij}^1(\omega), \hat{f}_{ij}^0(\omega)\right)\right)}{\sum_{i,j=1}^{p} \left( \operatorname{Var}(\hat{f}_{ij}^1(\omega) - \hat{f}_{ij}^0(\omega)) + \left|f_{ij}^1(\omega) - f_{ij}^0(\omega)\right|^2 \right)} \qquad (23)$$

*Proof.* The proof is found in A.1. □

### 3.6. Asymptotic behaviour of optimal shrinkage intensity

Now, we will examine the asymptotic behavior of the optimal shrinkage intensity (23). This will enable us to derive a data driven estimator in the following subsection. We first define the following parameters:

$$\pi(\omega) = \sum_{i,j=1}^{p} \pi_{ij}(\omega) \qquad (24)$$

$$\rho(\omega) = \sum_{i,j=1}^{p} \rho_{ij}(\omega) \qquad (25)$$

$$\gamma(\omega) = \sum_{i,j=1}^{p} \gamma_{ij}(\omega) \qquad (26)$$

where the subcomponents are defined, respectively, as:

$$\pi_{ij}(\omega) = \operatorname{AsyVar}\left(\sqrt{m_T} \hat{f}_{ij}^0(\omega)\right) = |f_{ij}(\omega)|^2 \qquad (27)$$

$$\rho_{ij}(\omega) = \operatorname{AsyCov}\left(\sqrt{m_T} \hat{f}_{ij}^1(\omega), \sqrt{m_T} \hat{f}_{ij}^0(\omega)\right) = \beta_i \beta_j f_{0i}(\omega) f_{j0}(\omega) \qquad (28)$$

$$\gamma_{ij}(\omega) = \left|f_{ij}^1(\omega) - f_{ij}^0(\omega)\right|^2 \tag{29}$$

using the notation

$$\text{AsyVar}(\cdot) := \lim_{T \to \infty} \text{Var}(\cdot)$$

and with weights $\beta_i$ defined in equation (7). Now, we can express $\zeta_T^*(\omega)$ as a function of (24) to (26) plus a remainder term which converges to zero sufficiently fast under the following additional assumption:

*Assumption* 3.4. The smoothing span $m_T$ is supposed to fulfill $m_T^2/T \to 0$ as $T \to \infty$.

This assumption 3.4 is made for the technical reason of proving the following theorem which gives now the exact expression of $\zeta_T^*(\omega)$:

*Theorem* 3.4. The optimal shrinkage intensity can be expressed as the following function of the parameters $\pi(\omega), \rho(\omega)$ and $\gamma(\omega)$:

$$\zeta_T^*(\omega) = \frac{1}{m_T} \frac{\pi(\omega) - 2\Re(\rho(\omega))}{\gamma(\omega)} + \mathrm{O}\left(T^{-1/2}\right) \tag{30}$$

*Proof.* The proof is found in A.2 □

This means that the optimal shrinkage intensity converges to zero at a rate of $1/m_T$. At the same time, it can be approximated by the parameters (24) to (26) with an error that vanishes, under assumption 3.4, with the faster rate of $T^{-1/2}$. This will allow us to derive a data driven estimator of the shrinkage intensity, and thus of $f_T^0(\omega)$, by plugging in estimators for (24) to (26) in (30).

### 3.7. Data driven estimation

The final step in deriving a data driven estimator of the spectrum that combines the averaged periodogram with a parsimonious, one-factor model based estimator, is to derive estimators for the parameters $\pi(\omega), \rho(\omega)$ and $\gamma(\omega)$. We will start by estimating $\pi(\omega)$. According to (24), $\pi_{ij}(\omega)$ is the asymptotic variance of the $i,j$th component of the averaged periodogram, scaled by the smoothing span $m_T$. The following lemma will give a consistent estimator:

*Lemma* 3.5. $\pi(\omega)$ is estimated consistently by

$$p(\omega) = \sum_{i,j=1}^{p} p_{ij}(\omega) \tag{31}$$

where

$$p_{ij}(\omega) = \frac{1}{m_T} \sum_{k=-(m_T-1)/2}^{(m_T-1)/2} |I_{ij}(\omega + \omega_k) - \hat{f}_{ij}^0(\omega)|^2 \tag{32}$$

i.e. $p_{ij}(\omega)$ is the standard estimator of the local variance of the $(i,j)$th component of the periodogram at frequency $\omega$.

*Proof.* The proof is given in A.3 □

The next step is to estimate $\rho(\omega)$. We will estimate its components and distinguish between the components on the diagonal and the components on the off-diagonal. As observed earlier, on the diagonal, $\hat{f}_T^1(\omega) = \hat{f}_T^0(\omega)$, thus we can use the estimator (32). On the off-diagonal, we can use the estimator given by the following lemma:

*Lemma* 3.6. For $i \neq j$, a consistent estimator of $\rho_{ij}(\omega)$ is given by

$$r_{ij}(\omega) = b_i b_j \hat{f}_{0i}^0(\omega) \hat{f}_{j0}^0(\omega) \tag{33}$$

*Proof.* The proof is given in A.4. □

The estimator of the last of the three parameters, $\gamma(\omega)$, is derived in a straightforward way:

*Lemma* 3.7. $\gamma(\omega)$ is estimated consistently by

$$g(\omega) = \sum_{i,j=1}^{p} g_{ij}(\omega) \tag{34}$$

where

$$g_{ij}(\omega) = \left| \hat{f}_{ij}^1(\omega) - \hat{f}_{ij}^0(\omega) \right|^2 \tag{35}$$

*Proof.* Both $\hat{f}_T^0(\omega)$ and $\hat{f}_T^1(\omega)$ are consistent estimators of $f_T^0(\omega)$ and $f_T^1(\omega)$, respectively. □

With the help of lemmata 3.5 to 3.7, we can now construct the data driven market shrinkage estimator of the spectrum, which is given by the following theorem:

*Theorem* 3.8. The estimator

$$\zeta_T(\omega) = \frac{1}{m_T} \frac{p(\omega) - 2\Re(r(\omega))}{g(\omega)} \tag{36}$$

is a consistent estimator of

$$\frac{1}{m_T} \frac{\pi(\omega) - 2\Re(\rho(\omega))}{\gamma(\omega)}.$$

*Proof.* This is implied by assumption 3.3 in conjunction with lemmata 3.5, 3.6 and 3.7. □

Thus, we have finally arrived at a shrinkage estimator that depends on the data only, not on background knowledge of the underlying process:

$$\hat{f}_T^+(\omega) = \frac{p(\omega) - 2\Re(r(\omega))}{m_T g(\omega)} \hat{f}_T^1(\omega) + \left(1 - \frac{p(\omega) - 2\Re(r(\omega))}{m_T g(\omega)}\right) \hat{f}_T^0(\omega) \tag{37}$$

We will refer to this estimator as to the DDMSE (data driven market shrinkage estimator). The following theorem gives the asymptotic behavior of the DDMSE $\hat{f}_T^+(\omega)$:

*Theorem* 3.9. Under assumptions 3.1 to 3.4, $\hat{f}_T^+(\omega)$ is a consistent estimator of the spectrum.

*Proof.* Asymptotically, the optimal shrinkage intensity $\zeta_T^*(\omega)$ vanishes according to theorem 3.4. According to theorem 3.8, $\zeta_T^*(\omega)$ is estimated consistently by $\zeta_T(\omega)$. This means that $\zeta_T(\omega)$ converges to zero, too, and thus that $\hat{f}_T^+(\omega)$ converges to the averaged periodogram. □

The performance of the DDMSE in practice will be examined by extensive Monte Carlo simulations in section 4.

## 4. Monte Carlo studies for the DDMSE

In this section, we will evaluate the performance of the data driven market shrinkage estimator in practice. For this, we will perform comprehensive Monte Carlo simulations. The DDMSE will have three benchmark estimators to compete with:

1. the averaged periodogram
2. the one-factor model that is the shrinkage target
3. a competing shrinkage estimator, referred to as *DDSSE*, that uses the identity matrix as the shrinkage target, see Böhm and von Sachs (2007)

In a setting where it is reasonable to use the DDMSE, it should outperform all three benchmarks. Such a setting can be characterized as the frequently encountered situation where one may fit a factor model to the data, but has no background knowledge on how many factors to actually choose. In a screeplot of the eigenvalues, one will typically encounter one or more prominent eigenvalues followed by a longer tail of small eigenvalues. The method we have developed will allow us to avoid the problem of model choice.

### 4.1. Setup

For the simulations, we have chosen to use a two-factor model as the true model. The first factor is an MA(2) process. Its spectrum has a peak at $\pi/2$. The second factor driving the process is a Gaussian white noise time series; its variance will be varied in a first simulation study, to examine the performance of the DDMSE on the 'scale' between almost one-factor model to true two-factor model. Figure 1 shows the spectrum of the two factors underlying the simulations. These two factors are then projected onto a 5-dimensional time series according to the following model:

$$X_t = \Upsilon f_t + \epsilon_t, \qquad \epsilon \sim \mathcal{N}(0, \Omega) \tag{38}$$

Here, $\Upsilon$ is a $5 \times 2$ weight matrix that was chosen at random initially, then fixed for this section. The initial random distribution for the components of $\Upsilon$ was uniform $\sim U([.3, 1])$, the components being chosen independently.

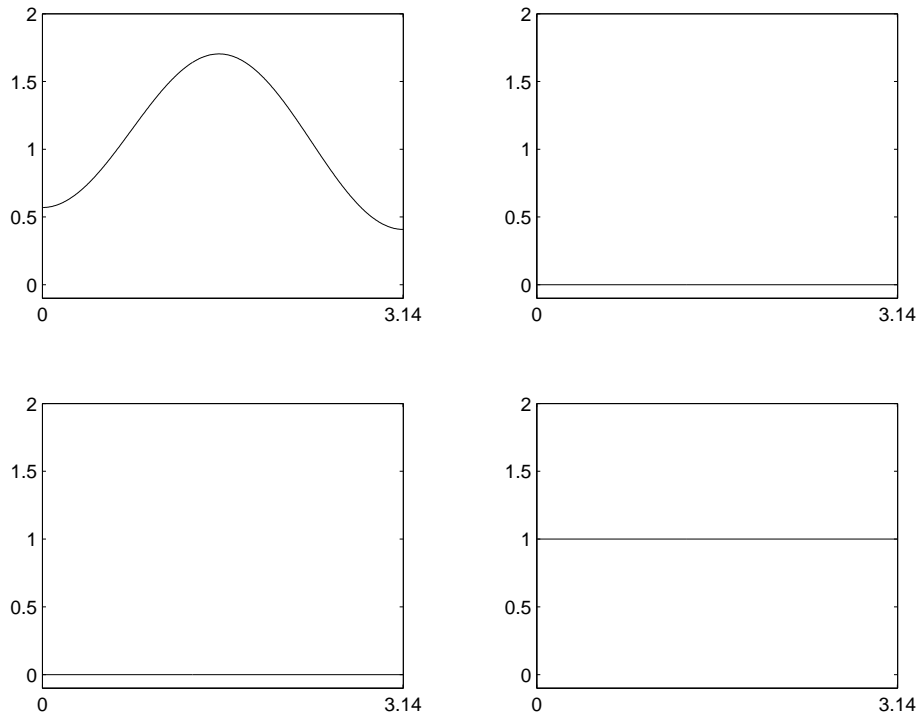

Fig 1: True spectrum of the two underlying factors. The imaginary parts are all zero.

$$\Upsilon = \begin{pmatrix} .5871 & .4510 \\ .5676 & .9691 \\ .4645 & .7268 \\ .8691 & .5511 \\ .5379 & .4754 \end{pmatrix}$$

The covariance matrix of the idiosyncratic components was obtained likewise: the off-diagonal components were set to zero, the diagonal components were simulated as iid uniform $\sim U([.2, .4])$ and then fixed as

$$\Omega = \text{diag}(.3213 \quad .3726 \quad .2646 \quad .4169 \quad .3257)$$

The market factor time series was obtained as the mean over dimension of the simulated data. All simulations presented in this section were repeated for new realizations of $\{\Upsilon, \text{Cov}(\Omega)\}$ without any major changes in the results, which is why we will omit these repetitive studies. A length of $T = 1,024$ was chosen for the time series in this section.

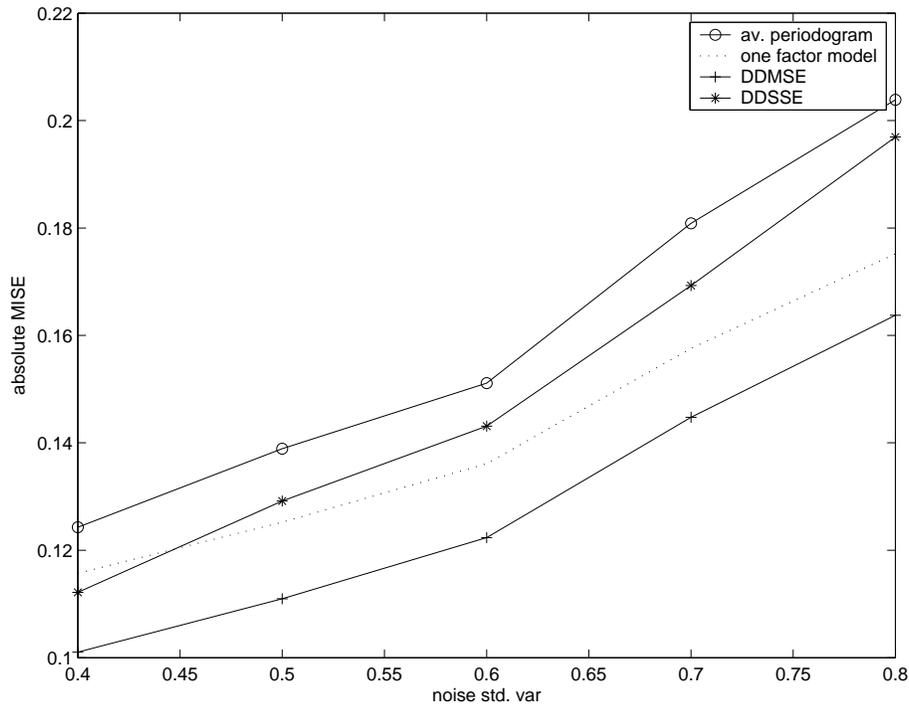

Fig 2: MISE of DDMSE, averaged periodogram, 1-factor-model and DDSSE for data from a 2 factor model. $T = 1,024$, smoothing span $m = 19$, different standard deviations of second factor. Based on $M = 1,500$ Monte Carlo runs. Confidence intervals are not printed as there is no intersection at 99% level for the solid curves.

### 4.2. Influence of the true model

The only formal prerequisite for the true model in order for the DDMSE to work is that its true spectrum is not that of a one-factor model such as the one specified in section 3.2. In this subsection, we will examine the influence of the 'distance' from a one-factor model. This is accomplished by using the two-factor model (38) to generate the data and systematically varying the standard deviation of the second, flat-spectrum factor. For small standard deviation, the data are very close to a one-factor model; as the standard deviation of the second factor increases, so does its influence. The results are given in figure 2. The effects we observe in the simulations study confirm our assumptions on the respective behavior of averaged periodogram, one-factor model, DDSSE and DDMSE. First of all, we remark that the DDMSE performs best for all choices of the white noise variance in the simulations. The averaged periodogram, upon which we want to improve, exhibits the worst performance. Not only is it outperformed by the DDSSE, which we would have expected based on the results of the preceding section, but also by the one-factor model. This shows that, in this context, the one-factor model is a useful model in itself, even although it is actually misspecified. It even outperforms the DDSSE for most choices of white noise variance. Overall, the MISE increases with the variance of the second factor, and the different estimators follow the MISE in a parallel shape.

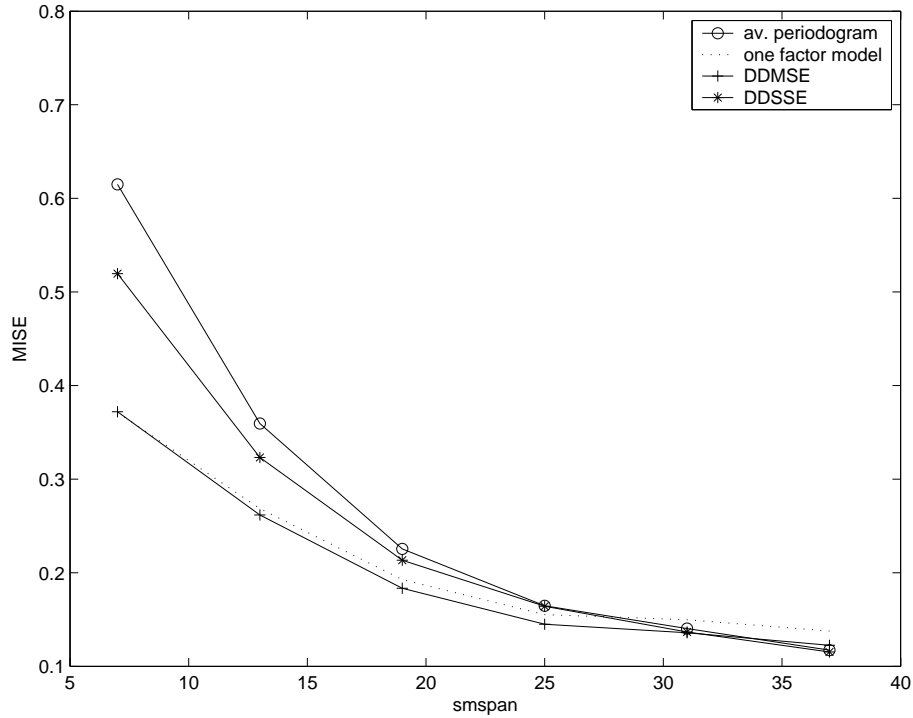

Fig 3: MISE of DDMSE, averaged periodogram, 1-factor-model and DDSSE for data from a 2 factor model. $T = 1,024$, different smoothing spans. Based on $M = 1,500$ Monte Carlo runs.

### *4.3. Influence of the smoothing span*

In the next Monte Carlo study, we have varied the smoothing span and examined its influence on the MISE. The results are given in figure 3. Not surprisingly, we observe that the overall MISE decreases as the sample size is increased for all three estimators. For small smoothing span, the averaged periodogram exhibits the worst performance. The DDSSE performs better than the averaged periodogram for small smoothing span, but is outperformed by the one-factor model and by the DDMSE. For the very small smoothing span $m = 7$, the DDMSE and the one-factor model have approximately the same MISE. Then, we have again the ranking averaged periodogram-DDSSE-one-factor model-DDMSE, as in the preceding subsection. Finally, for a comparatively large smoothing span of $m = 31$ or larger, the DDMSE, DDSSE and averaged periodogram seem to have approximately the same MISE. This is again not surprising, as for fixed dimension, both data driven estimators converge to the averaged periodogram. Moreover, for large smoothing span, the one-factor model performs worse than the averaged periodogram. This is, however, not due to a loss of performance of the one-factor model, which improves monotonously with $m$, but rather due to the faster improvement of the averaged periodogram in terms of MISE. Finally, the deterioration of the estimator based the one-factor model with respect to the averaged periodogram for large smoothing span does not make the DDMSE perform worse than the averaged periodogram. This can be explained by the fact that, for large $m$, the shrinkage intensity becomes negligibly small.

## 5. Conclusions

Our work deals with the concept of shrinkage in the *frequency domain* of multivariate time series. Similarly to our companion paper Böhm and von Sachs (2007), it uses a new, *localized* concept of shrinkage that allows for the development of estimators that simultaneously overcome the problem of numerical instability due to high dimensionality or collinearity and have lower quadratic risk. In contrast to the developments in the time domain of Ledoit and Wolf (2003), in the frequency domain of nonparametric estimation of the spectral density matrix by smoothing the periodogram matrix, all considerations have to be undertaken with respect to the (locally) effectively available sample size, which is governed by the smoothing parameter (and not the sample size alone). In Böhm and von Sachs (2007) asymptotic theory has been derived for the situation of shrinkage towards a multiple of the identity matrix where both the dimensionality $p = p_T$ and the smoothing span $m = m_T$ tend to infinity as the length of the time series $T \to \infty$. In this paper, we have contented ourselves to investigate the theoretical properties of our proposed estimator by classical asymptotics, noting that a transfer to the more complex situation of "Kolmogorov" or double asymptotics would be possible as well. However, with this work on structural shrinkage, we want to put emphasis onto a different aspect of shrinkage, perhaps driven by a more applied interest. Using the identity matrix as a shrinkage target is reasonable if there is little or no knowledge about the underlying multidimensional structure of the data. However, in many situations, in particular in economic applications, it is more rewarding to incorporate potentially available background knowledge on the underlying cross dimensional structure of the data into the shrinkage target. This opens up the way to designing 'custom made' shrinkage estimators that offer a new answer to problems of model choice. In a given setting where a class of parametric or semi-parametric estimators is eligible, and the order has to be chosen, instead of relying on criteria such as AIC or BIC, the minimum order model can be used as a target towards which to shrink. Instead of calling the method "shrinkage" we might as well describe it as *stretching*: a too parsimonious model is fitted and the estimate is then refined by adding the periodogram as a stretching target that has low bias and high variance.

In addition to showing that a MISE-optimal "oracle" shrinkage intensity can be consistently estimated from the data, we have shown by our Monte Carlo simulations, even for small sample size, the large gain in terms of $L_2$ risk of our estimator, in a situation of disposing additional structure, over the following competitors: the classical averaged periodogram, the "shrinkage to identity" estimator of Böhm and von Sachs (2007) and an estimator based on a fully parametric factor model. Simulations not reported here also demonstrate that shrinkage can be applied to tapered data; as tapering improves the rate of the bias without changing the rate of consistency, it is easy to transfer this to theory. For similar reasons, it is possible to replace the averaging of the periodogram by kernel smoothing.

An important field of application of our approach would be factor modelling of panels of economic time series data of comparatively high dimensionality. We recall that "high dimension" needs to be understood as high compared to the "effective sample size" $m_T$. Our achievements of this paper suggest that it could be possible to circumvent the problem of searching for an appropriate factor dimension - a problem still not satisfactorily solved in the literature, in particular for dynamic factor models. This latter application calls for a possible theoretical direction of future research: the generalization of our approach to a dynamic (and latent) factor model setting that allows for lag effects.

## 6. Acknowledgements

We would like to thank Christian Hafner and Johan Segers (UC Louvain) as well as Hernando Ombao (Brown University) for helpful discussions and an anonymous referee for his comments on ridge regression. Further, we acknowledge financial support from the IAP research network grant P 5/24 of the Belgian government (Belgian Science Policy) as well as from the contract "Projet d'Actions de Recherche Concertées" n 07/12-002 of the "Communauté française de Belgique", granted by the "Académie universitaire Louvain".

## Appendix A: Proofs

We will make frequent use of the abbreviations $\tilde{\omega}$, which means the Fourier frequency nearest $\omega$, and $\tilde{\omega}_k := \tilde{\omega} + \omega_k$.

### A.1. Proofs of equation (21) and of theorem 3.3

We begin by showing equation (21) which can be decomposed as follows:

$$R(z) = \mathrm{E}\left\|z\hat{f}_T^1(\omega) + (1-z)\hat{f}_T^0(\omega) - f_T^0(\omega)\right\|^2$$

$$\begin{aligned}
&= \sum_{i,j=1}^{p} \mathrm{E} \left| z \hat{f}_{ij}^1(\omega) + (1-z) \hat{f}_{ij}^0(\omega) - f_{ij}^0(\omega) \right|^2 \\
&= \sum_{i,j=1}^{p} \mathrm{Var} \left| z \hat{f}_{ij}^1(\omega) + (1-z) \hat{f}_{ij}^0(\omega) - f_{ij}^0(\omega) \right| \\
&\quad + \sum_{i,j=1}^{p} \left| \mathrm{E} \left( z \hat{f}_{ij}^1(\omega) + (1-z) \hat{f}_{ij}^0(\omega) - f_{ij}^0(\omega) \right) \right|^2 \\
&= \sum_{i,j=1}^{p} \mathrm{Var} \left| z \hat{f}_{ij}^1(\omega) + (1-z) \hat{f}_{ij}^0(\omega) - f_{ij}^0(\omega) \right| \\
&\quad + \sum_{i,j=1}^{p} \left| z f_{ij}^1(\omega) - z f_{ij}^0(\omega) \right|^2 \\
&= \sum_{i,j=1}^{p} \Big( z^2 \mathrm{Var}\, \hat{f}_{ij}^1(\omega) + (1-z)^2 \mathrm{Var}\, \hat{f}_{ij}^0(\omega) \\
&\qquad + z(1-z) \mathrm{Cov}\left( \hat{f}_{ij}^1(\omega), \hat{f}_{ij}^0(\omega) \right) \\
&\qquad + z(1-z) \mathrm{Cov}\left( \hat{f}_{ij}^0(\omega), \hat{f}_{ij}^1(\omega) \right) \\
&\qquad + z^2 \left| f_{ij}^1(\omega) - f_{ij}^0(\omega) \right|^2 \Big) \\
&= \sum_{i,j=1}^{p} \Big( z^2 \mathrm{Var}\, \hat{f}_{ij}^1(\omega) + (1-z)^2 \mathrm{Var}\, \hat{f}_{ij}^0(\omega) \\
&\qquad + 2z(1-z) \Re\left( \mathrm{Cov}\left( \hat{f}_{ij}^1(\omega), \hat{f}_{ij}^0(\omega) \right) \right) \\
&\qquad + z^2 \left| f_{ij}^1(\omega) - f_{ij}^0(\omega) \right|^2 \Big)
\end{aligned}$$

Then we want to derive the optimal shrinkage intensity $\zeta_T(\omega)$, which is the solution of the optimization problem (20). According to (22), any local extremum of the function $R(z)$ is a minimum. Thus, $\zeta_T^*(\omega)$ is the value obtained for $z$ by setting the first derivative equal to zero:

$$\begin{aligned}
&0 = R'(\zeta_T^*(\omega)) \\
\Leftrightarrow\quad &0 = 2 \sum_{i,j=1}^{p} \Big\{ \zeta_T^*(\omega) \mathrm{Var}\, \hat{f}_{ij}^1(\omega) - (1 - \zeta_T^*(\omega)) \mathrm{Var}\, \hat{f}_{ij}^0(\omega) \\
&\qquad + (1 - 2\zeta_T^*(\omega)) \Re\left( \mathrm{Cov}\left( \hat{f}_{ij}^1(\omega), \hat{f}_{ij}^0(\omega) \right) \right) \\
&\qquad + \zeta_T^*(\omega) \left| f_{ij}^1(\omega) - f_{ij}^0(\omega) \right|^2 \Big\} \\
\Leftrightarrow\quad &2\zeta_T^*(\omega) \sum_{i,j=1}^{p} \Big\{ \mathrm{Var}\, \hat{f}_{ij}^1(\omega) + \mathrm{Var}\, \hat{f}_{ij}^0(\omega)
\end{aligned}$$

$$-2\Re\left(\operatorname{Cov}\left(\hat{f}_{ij}^1(\omega),\hat{f}_{ij}^0(\omega)\right)\right)+\left|f_{ij}^1(\omega)-f_{ij}^0(\omega)\right|^2\Big\}$$

$$=2\sum_{i,j=1}^{p}\left\{\operatorname{Var}\hat{f}_{ij}^0(\omega)-2\Re\left(\operatorname{Cov}\left(\hat{f}_{ij}^1(\omega)-\hat{f}_{ij}^0(\omega)\right)\right)\right\}$$

$$\Leftrightarrow\quad 2\zeta_T^*(\omega)\sum_{i,j=1}^{p}\left\{\operatorname{Var}\left(\hat{f}_{ij}^1(\omega)-\hat{f}_{ij}^0(\omega)\right)+\left|f_{ij}^1(\omega)-f_{ij}^0(\omega)\right|^2\right\}$$

$$=2\sum_{i,j=1}^{p}\left\{\operatorname{Var}\hat{f}_{ij}^0(\omega)-2\Re\left(\operatorname{Cov}\left(\hat{f}_{ij}^1(\omega)-\hat{f}_{ij}^0(\omega)\right)\right)\right\}$$

$$\Rightarrow\quad \zeta_T^*(\omega)=\frac{\sum_{i,j=1}^{p}\left(\operatorname{Var}\hat{f}_{ij}^0(\omega)-2\Re\operatorname{Cov}\left(\hat{f}_{ij}^1(\omega),\hat{f}_{ij}^0(\omega)\right)\right)}{\sum_{i,j=1}^{p}\left(\operatorname{Var}\left(\hat{f}_{ij}^1(\omega)-\hat{f}_{ij}^0(\omega)\right)+\left|f_{ij}^1(\omega)-f_{ij}^0(\omega)\right|^2\right)}$$

$\square$

## A.2. Proof of theorem 3.4

Theorem 3.4 is proven using two technical lemmata which we will give immediately after the proof, which we give first:

If we multiply (23) by $m_T$, we obtain

$$m_T\zeta_T^*(\omega)=\frac{\sum_{i,j}\left(\operatorname{Var}\left(\sqrt{m_T}\hat{f}_{ij}^0(\omega)\right)-2\Re\left(\operatorname{Cov}\left(\sqrt{m_T}\hat{f}_{ij}^1(\omega),\sqrt{m_T}\hat{f}_{ij}^0(\omega)\right)\right)\right)}{\sum_{i,j=1}^{p}\left(\operatorname{Var}(\hat{f}_{ij}^1(\omega)-\hat{f}_{ij}^0(\omega))+\left|f_{ij}^1(\omega)-f_{ij}^0(\omega)\right|^2\right)} \tag{39}$$

$\hat{f}_{ij}^0(\omega)$ and $\hat{f}_{ij}^1(\omega)$ are consistent estimators of $f_{ij}^0(\omega)$ and $f_{ij}^1(\omega)$, respectively. This means that

$$\operatorname{Var}\left(\hat{f}_{ij}^1(\omega)-\hat{f}_{ij}^0(\omega)\right)=\operatorname{o}(1) \tag{40}$$

Using assumption 3.3 and (40), we obtain that the denominator of the right hand side of (39) is $\operatorname{O}(1)$. The numerator of the right hand side of (39) is $\pi(\omega)-2\Re(\rho(\omega))+\operatorname{O}\left(\frac{m_T}{\sqrt{T}}\right)$ according to lemmata A.1 and A.2. This yields

$$m_T\zeta_T^*(\omega)=\frac{\pi(\omega)+\rho(\omega)+\operatorname{O}\left(\frac{m_T}{\sqrt{T}}\right)}{\gamma(\omega)} \tag{41}$$

or, equivalently,

$$\zeta_T^*(\omega)=\frac{1}{m_T}\frac{\pi(\omega)+\rho(\omega)}{\gamma(\omega)}+\operatorname{O}\left(T^{-1/2}\right) \tag{42}$$

$\square$

*A.2.1. Lemmata needed for A.2 (proof of theorem 3.4)*

*Lemma* A.1.
$$\mathrm{Var}\left(\sqrt{m_T}\hat{f}^0_{ij}(\omega)\right) = \pi_{ij}(\omega) + \mathrm{O}\left(\frac{m_T}{T}\right) \qquad (43)$$

*Proof.*

$$\begin{aligned}
\mathrm{Var}(\sqrt{m_T}\hat{f}^0_{ij}(\omega)) &= \mathrm{Var}\left(\frac{1}{\sqrt{m_T}}\sum_{k=-(m_T-1)/2}^{(m_T-1)/2} I_{ij}(\tilde{\omega}_k)\right) \\
&= \frac{1}{m_T}\sum_{k=-(m_T-1)/2}^{(m_T-1)/2} \mathrm{Var}(I_{ij}(\tilde{\omega}_k)) \\
&\quad + \frac{1}{m_T}\sum_{\substack{k,l=-(m_T-1)/2 \\ k\neq l}}^{(m_T-1)/2} \mathrm{Cov}(I_{ij}(\tilde{\omega}_k)I_{ij}(\tilde{\omega}_l)) \\
&= |f_{ij}(\omega)|^2 + \mathrm{O}\left(\frac{m_T}{T}\right) + \frac{1}{m_T}\mathrm{O}\left(\frac{m_T^2}{T}\right) \\
&= |f_{ij}(\omega)|^2 + \mathrm{O}\left(\frac{m_T}{T}\right)
\end{aligned}$$

This proves equation (43) and yields that

$$\pi_{ij}(\omega) = |f_{ij}(\omega)|^2 \qquad (44)$$

□

*Lemma* A.2. For $i \neq j$,
$$\mathrm{Cov}\left(\sqrt{m_T}\hat{f}^1_{ij}(\omega), \sqrt{m_T}\hat{f}^0_{ij}(\omega)\right) = \rho_{ij}(\omega) + \mathrm{O}\left(\frac{m_T}{\sqrt{T}}\right). \qquad (45)$$

*Proof.* In the following estimate, we make use of (16), i.e. the convergence in probability of $b_i$ coming from equation (14),

$$b_i = \beta_i + \mathrm{O}_p\left(\frac{1}{\sqrt{T}}\right).$$

In order to control the error in replacing the random $b_i$ by their limiting $\beta_i$ we use Cauchy's inequality applied to all occuring remainder terms of the following or similar type

$$\mathrm{Cov}\left((b_i - \beta_i)\beta_j I_{00}(\tilde{\omega}_k), I_{ij}(\tilde{\omega}_l)\right).$$

With this we can derive that

$$\mathrm{Cov}\left(\sqrt{m_T}\hat{f}^1_{ij}(\omega), \sqrt{m_T}\hat{f}^0_{ij}(\omega)\right)$$
$$= \mathrm{Cov}\left(\frac{1}{\sqrt{m_T}}\sum_{k=-(m_T-1)/2}^{(m_T-1)/2} b_i b_j I_{00}(\tilde{\omega}_k), \frac{1}{\sqrt{m_T}}\sum_{k=-(m_T-1)/2}^{(m_T-1)/2} I_{ij}(\tilde{\omega}_k)\right)$$

$$= \frac{1}{m_T}\beta_i\beta_j \left\{ \sum_{k=-(m_T-1)/2}^{(m_T-1)/2} \text{Cov}\left(I_{00}(\tilde{\omega}_k), I_{ij}(\tilde{\omega}_k)\right) \right.$$

$$\left. + \sum_{\substack{k,l=-(m_T-1)/2 \\ k\neq l}}^{(m_T-1)/2} \text{Cov}\left(I_{00}(\tilde{\omega}_k), I_{ij}(\tilde{\omega}_l)\right) \right\} + \text{O}\left(\frac{m_T}{\sqrt{T}}\right)$$

$$= \frac{1}{m_T}\beta_i\beta_j \sum_{k=-(m_T-1)/2}^{(m_T-1)/2} \text{Cov}(I_{00}(\tilde{\omega}_k), I_{ij}(\tilde{\omega}_k)) + \text{O}\left(\frac{m_T}{T}\right) + \text{O}\left(\frac{m_T}{\sqrt{T}}\right) , \quad (46)$$

where we have used that $\text{Cov}\left(I_{00}(\tilde{\omega}_k), I_{ij}(\tilde{\omega}_l)\right) = \text{O}\left(\frac{1}{T}\right)$ for $k \neq \ell$. Showing this is parallel to treating the leading term, i.e. the covariance term in (46) using lemma A.3 and lemma A.4:

$$\text{Cov}(I_{00}(\omega), I_{ij}(\omega))$$
$$= \text{Cov}\left(d_0(\omega)\overline{d_0(\omega)}, d_i(\omega)\overline{d_j(\omega)}\right)$$
$$= \text{Cov}\left(d_0(\omega), d_i(\omega)\right) \text{Cov}\left(\overline{d_0(\omega)}, \overline{d_j(\omega)}\right)$$
$$\quad + \text{Cov}\left(d_0(\omega), \overline{d_j(\omega)}\right) \text{Cov}\left(\overline{d_0(\omega)}, d_i(\omega)\right)$$
$$= \text{E}\left(d_0(\omega)\, \overline{d_i(\omega)}\right) \text{E}\left(\overline{d_0(\omega)}\, d_j(\omega)\right) + \text{O}\left(\frac{1}{T}\right)\text{O}\left(\frac{1}{T}\right)$$
$$= f_{0i}(\omega)f_{j0}(\omega) + \text{O}\left(\frac{1}{T}\right) \quad (47)$$

Combining this with (46) and (47) yields thus

$$\text{Cov}\left(\sqrt{m_T}\hat{f}^1_{ij}(\omega), \sqrt{m_T}\hat{f}^0_{ij}(\omega)\right) = \beta_i\beta_j f_{0i}(\omega)f_{j0}(\omega) + \text{O}\left(\frac{m_T}{\sqrt{T}}\right) \quad (48)$$

which proves (45) and at the same time yields

$$\rho_{ij}(\omega) = \beta_i\beta_j f_{0i}(\omega)f_{j0}(\omega) . \quad (49)$$

□

### A.3. Proof of lemma 3.5

*Proof.* According to (Brockwell and Davis, 1987, theorem 10.3.2), we have

$$\text{Var}\, I_{ij}(\tilde{\omega}_k) = |f_{ij}(\tilde{\omega}_k)|^2 + \text{O}\left(1/\sqrt{T}\right) \quad (50)$$

and

$$\text{Cov}(I_{ij}(\tilde{\omega}_k), I_{ij}(\tilde{\omega}_l)) = \text{O}\left(1/T\right) \quad (51)$$

for $0 < \tilde{\omega}_k \neq \tilde{\omega}_l < \pi$. Furthermore, $\hat{f}_{ij}^0(\omega) = \mathrm{E}\, I_{ij}(\tilde{\omega}_k) + \mathrm{O}\left(m_T/T\right) + \mathrm{O}_p\left(1/m_T\right) = o_p(1)$ for all $\tilde{\omega}_k$ in the span of $m_T$.

This allows us to write

$$p_{ij}(\omega)$$
$$= \frac{1}{m_T} \sum_{k=-(m_T-1)/2}^{(m_T-1)/2} \left| I_{ij}(\tilde{\omega}_k) - \hat{f}_{ij}^0(\omega) \right|^2$$
$$= \frac{1}{m_T} \sum_{k=-(m_T-1)/2}^{(m_T-1)/2} \left| I_{ij}(\tilde{\omega}_k) - \mathrm{E}\, I_{ij}(\tilde{\omega}_k) + o_p(1) \right|^2$$
$$= \frac{1}{m_T} \sum_{k=-(m_T-1)/2}^{(m_T-1)/2} \left| I_{ij}(\tilde{\omega}_k) - \mathrm{E}\, I_{ij}(\tilde{\omega}_k) \right|^2 + o_p(1) ,$$

having used that $\left| I_{ij}(\tilde{\omega}_k) - \mathrm{E}\, I_{ij}(\tilde{\omega}_k) \right|^2$ is $\mathrm{O}_p(1)$. It remains to show that

$$\frac{1}{m_T} \sum_{k=-(m_T-1)/2}^{(m_T-1)/2} \left| I_{ij}(\tilde{\omega}_k) - \mathrm{E}\, I_{ij}(\tilde{\omega}_k) \right|^2$$

converges to $|f_{ij}(\omega)|^2 = \pi_{ij}(\omega)$ in probability.

We observe that (50) allows to control convergence of the mean, whereas we can control the variance by borrowing strength from a proof of a CLT for $\frac{1}{m_T} \sum_{k=-(m_T-1)/2}^{(m_T-1)/2} I_{ij}(\tilde{\omega}_k)$. One technique, frequently used and to be found, e.g., in Brillinger (1975), proof of Theorem 7.4.4., is to show that the cumulants of higher order than 2 of $\sqrt{m_T} \frac{1}{m_T} \sum_{k=-(m_T-1)/2}^{(m_T-1)/2} I_{ij}(\tilde{\omega}_k)$ tend to zero, i.e. in particuler the cumulants of order 4. But this includes in particular that

$$\left(\frac{1}{m_T}\right)^2 \sum_k \sum_\ell \mathrm{Cov}\left( I_{ij}^2(\tilde{\omega}_k), I_{ij}^2(\tilde{\omega}_l) \right) \to 0 ,$$

which is what is needed here. (An explicit calculation of this covariance would also be possible by application of Brillinger (1975), Theorem 4.3.1, using the fact that all but the second-order cumulants of order one up to eight of the occurring Gaussian mean zero $d_i(\omega)$ are identical zero, and that the second-order cumulants are products of expressions of the form $\mathrm{Cov}\left(d_i(\tilde{\omega}_k), d_j(\tilde{\omega}_l)\right)$ which tend to zero for $k \neq \ell$.) □

### *A.4. Proof of lemma 3.6*

*Proof.* $b_i$, $b_j$, $\hat{f}_{0i}^0(\omega)$ and $\hat{f}_{j0}^0(\omega)$ are consistent estimators of $\beta_i, \beta_j, f_{0i}(\omega)$ and $f_{j0}(\omega)$. This yields, in conjunction with (49), the result. □

### *A.5. Additional lemmata*

*Lemma* A.3. Let $(X_1, X_2, X_3, X_4)$ be a 4-variate normal random variable. Then we have

$$\begin{aligned}&\mathrm{Cov}(X_1 X_2, X_3 X_4)\\ =\ &\mathrm{Cov}(X_1, X_3)\mathrm{Cov}(X_2, X_4) + \mathrm{Cov}(X_1, X_4)\mathrm{Cov}(X_2, X_3)\end{aligned}$$

*Proof.* The proof is found in Brillinger (1975, p. 21). □

*Lemma* A.4. For $i \neq j$, we have that

$$\mathrm{E}\, d_i(\omega_1) d_j(\omega_2) = \mathrm{O}\left(\frac{1}{T}\right) \qquad (52)$$

where the null convergence is uniform in $\{\omega_1, \omega_2\} \in (0, 2\pi] \times (0, 2\pi]$.

*Proof.* The proof of this can be found in Shumway and Stoffer (2000, p. 275ff). □